\documentclass[11pt]{amsart}
\usepackage[a4paper, margin=2.54cm]{geometry}

\usepackage[english]{babel} 
\usepackage{microtype}  
\usepackage{amsmath, amssymb, amsthm}  
\usepackage{mathtools}  
\usepackage{centernot}  
\usepackage{todonotes}  
\usepackage{enumitem}   
\usepackage{graphicx}   
\usepackage{xcolor}     
\usepackage{float}      

\usepackage[english=american,autostyle]{csquotes}  
\MakeOuterQuote{"}

\usepackage[
    backend=biber,
    style=alphabetic,
    giveninits=false, 
    isbn=false,
    doi=false,
    url=true,
    maxbibnames=10,
]{biblatex}

\DeclareFieldFormat[article,inbook,incollection,inproceedings,patent,thesis,unpublished]{title}{#1}
\DeclareFieldFormat[inbook, article,inproceedings,incollection,unpublished,thesis]{title}{\mkbibitalic{#1}}
\DeclareFieldFormat[book]{title}{#1}
\DeclareFieldFormat[inproceedings]{booktitle}{\mkbibitalic{#1}}
\DeclareFieldFormat[inbook, incollection]{booktitle}{#1}

\DefineBibliographyStrings{english}{
  byeditor = {edited by},
  editor   = {editor},
  editors  = {editors},
  chapter = {chapter}
}
\DeclareFieldFormat[article]{volume}{\textbf{#1}}
\DeclareFieldFormat[article]{number}{(#1)}
\renewbibmacro*{journal+issuetitle}{%
  \usebibmacro{journal}%
  \setunit*{\addspace}%
  \printfield{volume}%
  \printfield{number}%
  \setunit{\addcomma\space}%
  \usebibmacro{date}%
  \setunit{\addcomma\space}%
  \usebibmacro{issue}%
  \newunit}
\DeclareFieldFormat{url}{\url{#1}}

\theoremstyle{plain}
\newtheorem{theorem}{Theorem}  

\newtheorem{lemma}[theorem]{Lemma}

\newtheorem{observation}[theorem]{Observation}
\newtheorem{corollary}[theorem]{Corollary}

\theoremstyle{definition}
\newtheorem{definition}[theorem]{Definition}

\newtheorem*{example*}{Example}

\theoremstyle{remark}

\usepackage[unicode]{hyperref} 
\hypersetup{pdftitle=An Elementary Proof of Stirling's Formula}
\hypersetup{pdfauthor=Jakub Smolík}

\title{An Elementary Proof of Stirling's Formula}
\author{Jakub Smolík}
\date{October 2023}

\begin{document}
\maketitle

We give a short elementary proof of Stirling's formula $n! \sim \sqrt{2\pi n}\,n^n e^{-n}$. I learned about this proof from Dominic Beck, and it also appears in~\cite{little2015real}. The proof relies on the Wallis product which itself can be derived only using elementary algebra~\cite{wastlund2007elementary}.

\begin{lemma}
    The sequence $\displaystyle a_n \coloneqq \frac{n!}{\sqrt{n}\,n^ne^{-n}}$ 
    converges to a positive number.
\end{lemma}
\begin{proof}
    First we will show, that the sequence $(a_n)$ is decreasing and then prove, that it is bounded from below by a positive number. To begin, note that
    \[
        \frac{a_n}{a_{n+1}} = \frac{n!}{\sqrt{n}\,n^ne^{-n}} \cdot \frac{\sqrt{n+1}\,(n+1)^{n+1}\,e^{-n-1}}{(n+1)n!} =
        \frac{\sqrt{n+1}}{\sqrt{n}} \cdot \frac{(n+1)^n}{n^n e} =
        \frac{1}{e} \biggl(\frac{n+1}{n} \biggr)^\frac{2n+1}{2}.
    \]
    Now define $b_n \coloneqq \ln(a_n)$. Then
    \[
        b_n - b_{n+1} = \ln\biggl( \frac{a_n}{a_{n+1}} \biggr) = 
        \frac{2n+1}{2} \ln\biggl( \frac{n+1}{n} \biggr) - 1.
    \]
    Next we introduce a new variable, $k$, such that $\frac{n+1}{n} = \frac{1+k}{1-k}$. This choice of $k$ proves useful, as it allows us to utilize a Taylor series expansion. To satisfy this condition we set $k \coloneqq \frac{1}{2n+1}$, leading to the following:

    \[
        b_n - b_{n+1} = \frac{2n+1}{2} \ln\biggl( \frac{n+1}{n} \biggr) - 1 = \frac{1}{2k}\ln\biggl( \frac{1+k}{1-k} \biggr) - 1.
    \]
    Using a Taylor series expansion, we get        
    \begin{equation*}
    \begin{split}
        \ln\biggl( \frac{1+k}{1-k} \biggr) &= \ln(1+k) - \ln(1-k) = \\
        &= \biggl(k - \frac{k^2}{2} + \frac{k^3}{3} - \frac{k^4}{4} + \cdots \biggr) - \biggl(-k - \frac{k^2}{2} - \frac{k^3}{3} - \frac{k^4}{4} - \cdots \biggr) = \\
        &= 2\biggl( k + \frac{k^3}{3} + \frac{k^5}{5} \cdots \biggr) = 
        2\sum_{i=0}^\infty \frac{k^{2i+1}}{2i+1}.
    \end{split}
    \end{equation*}
    Now we have 
    \[
        b_n - b_{n+1} = \frac{1}{2k}\ln\biggl(\frac{1+k}{1-k}\biggr) - 1 = \sum_{i=0}^\infty \frac{k^{2i}}{2i+1} - 1 = \sum_{i=1}^\infty \frac{k^{2i}}{2i+1} > 0.
    \]
    Therefore, the sequence $(b_n)$ is decreasing; the natural logarithm is a monotonic function and so $(a_n)$ is decreasing as well. In order to show that it is bounded from bellow, we resume the calculation, noting that $0 < k < 1$:
    \begin{allowdisplaybreaks}
    \begin{equation*}
    \begin{split}
        b_n - b_{n+1} &= \sum_{i=1}^\infty \frac{k^{2i}}{2i+1} < \sum_{i=1}^\infty k^{2i} = k^2\sum_{i=1}^\infty k^{2i-2} = k^2\sum_{i=0}^\infty k^{2i} = \\
        &= \frac{k^2}{1-k^2} = \frac{1}{(2n+1)^2\Bigl( 1 - \frac{1}{(2n+1)^2} \Bigr)} = 
        \frac{1}{(2n+1)^2 - 1} = \\
        &= \frac{1}{2n(2n+2)} = \frac{1}{4n(n+1)} = \frac{1}{4n} - \frac{1}{4(n+1)}.
    \end{split}
    \end{equation*}
    \end{allowdisplaybreaks}
    Hence 
    \[
        b_n - \frac{1}{4n} < b_{n+1} - \frac{1}{4(n+1)}.
    \]
    We see that the sequence $\bigl(b_n - \frac{1}{4n}\bigr)$ is increasing, therefore 
    \[
        b_n > b_n - \frac{1}{4n} > b_1 - \frac{1}{4} = \frac{3}{4} \implies a_n > e^{0.75},
    \]
    bounding $(a_n)$ from bellow. This completes the proof.
\end{proof}

 We have shown, that $n!$ grows up to a constant multiple as does $\sqrt{n}\,n^ne^{-n}$. We will need the following lemma to find this constant.

\begin{definition}
    Define $(2n)!! \coloneqq 2\cdot4\cdot6\cdots(2n)$ and $(2n-1)!! \coloneqq 1\cdot3\cdot5\cdots(2n-1)$. 
\end{definition}
\begin{observation}
    It holds, that $(2n)!!(2n-1)!! = (2n)!$ and $(2n)!! = 2^n\,n!$. 
\end{observation}

\begin{lemma}
    $\displaystyle \lim_{n\to\infty} \frac{4^n\,n!^2}{\sqrt{n}\,(2n)!} = \sqrt{\pi}.$
\end{lemma}
\begin{proof}
    Rewrite the famous Wallis product $\displaystyle \frac{2}{1}\cdot \frac{2}{3}\cdot \frac{4}{3}\cdot \frac{4}{5}\cdot \frac{6}{5}\cdot \frac{6}{7} \cdots = \frac{\pi}{2}$ as a limit:
    \[
    \lim_{n\to\infty} \frac{(2n)!!^2}{(2n-1)!!^2(2n+1)} = \frac{\pi}{2}.
    \]
    By expanding the fraction and simplifying the double factorials, we obtain
    \[
    \lim_{n\to\infty} \frac{(2n)!!^2 (2n)!!^2}{(2n)!!^2 (2n-1)!!^2(2n+1)} = \lim_{n\to\infty} \frac{2^{4n}\,n!^4}{(2n)!^2(2n+1)} = \frac{\pi}{2}.
    \]
    After taking the square root, we get 
    \[
    \lim_{n\to\infty} \frac{4^n\,n!^2}{(2n)!\sqrt{2n+1}} = \sqrt{\frac{\pi}{2}}
    \implies \lim_{n\to\infty} \frac{4^n\,n!^2 \sqrt{2}}{(2n)!\sqrt{2n+1}} = \lim_{n\to\infty} \frac{4^n\,n!^2}{(2n)!\sqrt{n}} = \sqrt{\pi}.\qedhere
    \]
\end{proof}

\begin{theorem}
    $\displaystyle A \coloneqq \lim_{n\to\infty} a_n = \lim_{n\to\infty} \frac{n!}{\sqrt{n}\,n^ne^{-n}} = \sqrt{2\pi}.$
\end{theorem}
\begin{proof}
    We start with the limit derived in the preceding lemma. We have
    \begin{equation*}
    \begin{split}
         \sqrt{\pi} &= \lim_{n\to\infty} \frac{4^n\,n!^2}{\sqrt{n}\,(2n)!} = \lim_{n\to\infty} \frac
        {4^n \biggl(\displaystyle\frac{n!}{\sqrt{n}\,n^ne^{-n}} \biggr)^2 \bigl( \sqrt{n}\,n^ne^{-n} \bigr)^2}
        {\sqrt{n}\, \displaystyle\frac{(2n)!}{\sqrt{2n}\,(2n)^{2n}e^{-2n}} \sqrt{2n}\,(2n)^{2n}e^{-2n}} =\\
        &= \lim_{n\to\infty} \frac{4^n\,A^2\,n\,n^{2n}e^{-2n}}{\sqrt{n}\,A\,\sqrt{2n}\,(2n)^{2n}e^{-2n}} = \lim_{n\to\infty}\frac{A}{\sqrt{2}}.
    \end{split}
    \end{equation*}
    Therefore $A = \sqrt{2\pi}.$ 
\end{proof}

\begin{corollary} Thus, we obtain Stirling's formula $\displaystyle n! \sim \sqrt{2\pi n}\Bigl(\frac{n}{e}\Bigr)^n.$ 
\end{corollary}

\emergencystretch=1em
\printbibliography
\end{document}